\documentclass[a4paper,12pt]{article}
\usepackage[english]{babel}
\usepackage[T2A]{fontenc}
\usepackage{amsthm}
\usepackage[tbtags]{amsmath}
\usepackage{amsfonts,amssymb}
\begin{document}

\newtheorem{theorem}{Theorem}
\newtheorem{lemma}{Lemma}
\newtheorem{proposition}{Proposition}
\newtheorem{Cor}{Corollary}

\begin{center}{\large\bf Seminoetherian Modules over Non-Primitive HNP rings}
\end{center}
\begin{center}
Askar Tuganbaev\footnote{National Research University MPEI, Moscow 111250, Krasnokazarmennaya street, 14-1, Moscow; Lomonosov Moscow State University, Russia 119991, Moscow, Lenin Hills, 1; email: tuganbaev@gmail.com .}
\end{center}

\textbf{Abstract.} We study the structure of seminoetherian modules. Seminoetherian modules over non-primitive hereditary noetherian prime rings are completely described.


The study is supported by grant of Russian Science Foundation (=RSF).

\textbf{MSC2020 database 16D10, 16D70}

\section{Introduction}\label{section1}

$\,$

We only consider associative unital non-zero rings and unitary modules. The words of type <<a noetherian ring $A$>> mean that the both modules $A_A$ and $_AA$ are noetherian. The word <<$A$-module>> usually means <<a right $A$-module>>. 
All notions, which are not defined in the paper and used in it, are standard; e.g., see
\cite{CarE56}, \cite{Fai76}, \cite{Goo76}, \cite{GooW89}, \cite{Kas82}, \cite{Lam99}, \cite{Tug98}, \cite{Wis91}.

A module is called a \textbf{max module} if every its non-zero submodule has a maximal submodule; for example, max modules were considered in \cite{Sho74}, \cite{Cam75}, \cite{Cam75b}, \cite{Fai95}, \cite{Tug02}. A module is said to be \textbf{seminoetherian} is all its factor modules are max modules, i.e., if all its non-zero subfactors\footnote{Subfactors are submodules of factor modules.} have maximal submodules; see \cite[Chapter 1; 1.35, 1.36]{Tug98}. Every noetherian module is seminoetherian. All semisimple modules are seminoetherian but they are noetherian only if they are finitely generated.

It is known that the theory of modules over non-primitive hereditary noetherian prime rings is in many ways reminiscent of the theory of abelian groups ; see \cite{Sin74}, \cite{Sin75}, \cite{Sin76}, \cite{Sin78}, \cite{Sin79}, \cite{SinT79}.
The main result of this paper is Theorem 1.11\footnote{Apparently, this result is new even for abelian groups.}; the corresponding definitions and notations are given in Introduction. 

\textbf{1.1. Theorem.} Let $A$ be a non-primitive hereditary noetherian prime ring, $M$ be a right $A$-module, and let $T$ be the largest singular submodule of the module $M$. The following conditions are equivalent.
\begin{enumerate}
\item[\textbf{1)}]
$M$ is a seminoetherian module. 
\item[\textbf{2)}]
Every primary component of the module $T$ is a direct sum of cyclic uniserial modules, the composition lengths of which are limited in total, $M/T$ is a non-singular finite-dimensional (in the sense of Goldie) module and there exists a submodule $X$ of the module $M$ such that $T\subseteq X$, $X/T$ is a noetherian projective essential submodule of $M/T$ and every primary component of the module $M/X$ is a direct sum of cyclic uniserial modules, the composition lengths of which are limited in total.
\end{enumerate}

The proof of Theorem 1.1 is divided to several assertions from Sections 2, 3 and 4 of this paper.

A submodule $X$ of the module $M$ is said to be \textbf{maximal} if the factor module $M/X$ is simple, i.e.,. $M/X$ is the unique non-zero submodule of $M/X$. A ring $A$ is said to be \textbf{right primitive} (resp., \textbf{left primitive}) if there exists a faithful simple right (resp., left) $A$-module. A ring is said to be \textbf{non-primitive} if it is not right or left primitive. A module $M_A$ is said to be \textbf{projective} if it satisfies the following equivalent conditions: \textbf{1)} for any module $X_A$, every epimorphism $h\colon X\to \overline{X}$, and an arbitrary homomorphism $\overline{f}\colon M\to \overline{X}$, there exists a homomorphism $f\colon M\to X$ with $\overline{f}=hf$; \textbf{2)} the module $M$ is isomorphic to a direct summand of some direct sum of copies of the module $A_A$. A module is said to be \textbf{hereditary} if all its submodules are projective. A module $M$ is said to be \textbf{noetherian} if all its submodules are finitely generated (equivalently, $M$ is a module with maximum condition on ascending chains of submodules). A module $M$ is said to be \textbf{finite-dimensional} (in the sense of Goldie) if $M$ does not have submodules which are infinite direct sums of non-zero modules. A submodule $X$ of the module $M$ is said to be \textbf{essential} if $X\cap Y\ne 0$ for every non-zero submodule $Y$ of $M$. A ring $A$ is said to be \textbf{right bounded} (resp., \textbf{left bounded}) if every its essential right (resp., \textbf{left}) ideal contains a non-zero ideal of the ring $A$. A ring is said to be \textbf{prime} if the product of any two its non-zero ideals is not equal to the zero. A ring $A$ is said to be \textbf{semiprime} if $A$ does not have non-zero nilpotent ideals. A right finite-dimensional ring with maximum condition on right annihilators is called a \textbf{right Goldie ring}. If $N$ is a subset of the right (resp., left) $A$-module $M$, then $r(N)$ (resp., $\ell(N)$) denotes the right (resp., left) annihilator of the set $N$ in $A$.

Without further specification, we will use the known facts listed in Remarks 1.2--1.8 below.

\textbf{1.2. Remark.} Hereditary noetherian prime rings also are called \textbf{HNP rings}. In \cite{Len73}, it is proved that any HNP ring $A$ is (right and left) primitive or (right and left) bounded; in addition, if $A$ is a primitive bounded ring, then $A$ is a simple artinian ring. Therefore, non-primitive HNP rings coincide with non-artinian bounded HNP rings.

\textbf{1.3. Remark.} We denote by $\text{Sing }M$ the set of all elements $m$ of the right (resp., left) $A$-module $M$ such that $r(m)$ (resp., $\ell(m)$) is an essential right (resp., left) ideal of the ring $A$. It is well known that $\text{Sing }M$ is a fully invariant submodule of the module $M$; it is called the \textbf{singular submodule} of the module $M$. If $\text{Sing }M=M$ (resp., $\text{Sing }M=0$), then the module $M$ is said to be \textbf{singular} (resp., \textbf{non-singular}). We denote by $T(M)$ the set of all elements $m$ of the module $M_A$ such that $r(m)$ contains a regular\footnote{An element $a$ of the ring $A$ is said to be \textbf{regular} if $r(a)=\ell(a)=0$.} element of the ring $A$; this set is called the \textbf{torsion part} of the module $M$. A module $M$ is called a \textbf{torsion module} (resp., a \textbf{torsion-free module}) if $T(M)=M$ (resp., $T(M)=0$). 

All assertions of the following Remark 1.4 are well known; see, e.g., \cite[5.9,~5.10,~6.14,~6.10(a)]{GooW89}.

\textbf{1.4. Remark.} Let $A$ be a semiprime right Goldie ring. Then $A$ has the semisimple artinian right classical ring of fractions $Q$ and non-zero injective non-singular indecomposable right $A$-modules coincide (up to an $A$-module isomorphism) with minimal right ideals of the semisimple artinian ring $Q$ (if the ring $A$ is prime, then $Q$ is a simple ring which is isomorphic to a matrix ring over a division ring). 
The set of all essential right ideals of the ring $A$ coincides with the set of all right ideals of the ring $A$ containing regular elements. Consequently, the class of all singular (resp., non-singular) right $A$-modules coincides with the class of all torsion right $A$-modules (resp., torsion-free $A$-modules). All essential extensions of singular (resp., non-singular) right $A$-modules are singular (resp., non-singular) modules. For any module $M_A$, the module $M/\text{Sing}(M)$ is non-singular. 

Remark 1.4 implies Remark 1.5.

\textbf{1.5. Remark.} If $A$ is a semiprime right Goldie ring, then every non-torsion right $A$-module contains a non-zero non-singular submodule.

\textbf{1.6.} Let $A$ be a semiprime right Goldie ring, and let $Q$ be semisimple artinian right classical ring of fractions of the ring $A$ . An ideal $B$ of the ring $A$ is said to be {\bf invertible} if there exists a subbimodule $B^{-1}$ of the bimodule ${}_AQ_A$ such that $BB^{-1}=B^{-1}B=A$. Maximal elements of the set of all proper invertible ideals of the ring $A$ are called {\bf maximal invertible} ideals. The set of all maximal invertible ideals of the ring $A$ is denoted by ${\mathcal P}(A)$. If $M$ is an $A$-module and $P\in {\mathcal P}(A)$, then the submodule $\{m\in M\;\vert \;mP^n=0, ~n=1,2, \ldots \}$ is called the {\bf $P$-primary component} of the module $M$; it is denoted by $M(P)$. If $M=M(P)$ for some $P\in {\mathcal P}(A)$, then $M$ is called a {\bf primary} module or a {\bf $P$-primary} module.

\textbf{1.7. Remark.} Let $A$ be a non-primitive HNP ring, $M$ be a singular $A$-module, and let $\{M(P_i)\}_{i\in I}$ be the set of all primary components of the module $M$. Many properties of the singular module $M$ and of its primary components are well known; e.g., see \cite{Sin74}, \cite{Sin75}, \cite{Sin76}, \cite{Sin78}, \cite{Sin79}, \cite{SinT79}. These properties are similar to properties of primary components of torsion abelian groups. For example, every uniform $A$-module and, in particular, every indecomposable injective $A$-module is $P$-primary for some maximal invertible ideal $P$ of the ring $A$, $M=\oplus_{i\in I}M(P_i)$, all primary components $M(P_i)$ are fully invariant in $M$, for any submodule $X$ of the module $M$, we have $X=\oplus_{i\in I}X\cap M(P_i)$, where $X\cap M(P_i)$ are primary components of the module $X$, $M/X=\oplus_{i\in I}M/X\cap M(P_i$, and so on.

A module $M$ is said to be \textbf{uniform} if any two its non-zero submodules have the non-zero intersection, i.e., if any non-zero submodule of $M$ is essential. A module $M$ is said to be \textbf{uniserial} if the lattice of all its submodules is a chain, i.e., any two submodules of $M$ are comparable with respect to inclusion (equivalently, any two cyclic submodules of $M$ are comparable with respect to inclusion). Every uniserial module is uniform and every uniform module is indecomposable. A module $M$ is said to be \textbf{infinite-dimensional} if $M$ contains a submodule which is an infinite direct sum of non-zero modules. A module $M$ is said to be \textbf{injective} if $M$ is a direct summand of any module containing $M$. A right (resp., left) module $M$ is said to be \textbf{faithful} if $r(M)=0$ (resp., $\ell(M)=0$). 
 
\section{Properties of Seminoetherian Modules}\label{sec2}

\textbf{2.1. Proposition.} For a module $M$, the following conditions are equivalent.
\begin{enumerate}
\item[\textbf{1)}]
$M$ is a seminoetherian module. 
\item[\textbf{2)}]
$X$ and $M/X$ are seminoetherian modules for every submodule $X$ of the module $M$.
\item[\textbf{3)}]
$X$ and $M/X$ are seminoetherian modules for some proper submodule $X$ of the module $M$.
\end{enumerate}

\textbf{Proof.} Implications 1)\,$\Rightarrow$\,2)\,$\Rightarrow$\,3) are directly verified.

3)\,$\Rightarrow$\,1). Let $Y$ be a proper submodule of $M$ and let $Y/Z$ be a non-zero subfactor of the module $M$. We have to prove that $Y/Z$ is contained in some maximal submodule $N/Z$ of the module $M/Z$. Three cases are only possible: \textbf{1}, \textbf{2}, and \textbf{3}.

\textbf{Case 1:} $Y$ is a submodule of $X$. Since the non-zero seminoetherian module $M/X$ has some maximal submodule $N/X$, the module $X$ is contained in the maximal submodule $N$ of the module $M$. Therefore, the module $X/Z$ is contained in the maximal submodule $N/Z$ of the module $M/Z$. Then $Y/Z\subseteq N/Z$ .

\textbf{Case 2:} $X$ is a proper submodule of $Y\subsetneq M$. The proper submodule $Y/X$ of the seminoetherian module $M/X$ is contained in some maximal submodule $N/X$ of the module $M/X$. Therefore, the module $Y$ is contained in the maximal submodule $N$ of the module $M$.

\textbf{Case 3:} The modules $X$ and $Y$ are not comparable with respect to inclusion, whence $X$ and $Y$ are proper submodules of $X+Y$. The following two cases \textbf{3$_a$} and \textbf{3$_b$} are only possible.

 \textbf{Case 3$_a$:} $X+Y$ is a proper submodule of $M$. Since the non-zero seminoetherian module $M/(X+Y)$ has some maximal submodule $N/(X+Y)$, the module $X+Y$ is contained in the maximal submodule $N$ of the module $M$. Then $Y$ is contained in the maximal submodule $N$ of the module $M$. Therefore, the module $Y/Z$ is contained in the maximal submodule $N/Z$ of the module $M/Z$.
 
 \textbf{Case 3$_b$:} $X+Y=M$. Then the non-zero module $M/Y= (X+Y)/Y$ is seminoetherian, since it is isomorphic to the non-zero seminoetherian module $X/(X\cap Y)$. Therefore, the module $M/Y$ has the maximal submodule $N/Y$. Then $Y/Z$ is contained in the maximal submodule $N/Z$ of the module $M/Z$.~\hfill$\square$

\textbf{2.2. Corollary.} Let a module $M$ contain a finite chain of submodules
$$
0=X_0\subseteq X_1\subseteq\dots\subseteq X_n\subseteq X_{n+1}=M
$$
such that the modules $M_{k+1}/M_k$ are seminoetherian for $k=0,\dots,n$. Then the module $M$ is seminoetherian.

Corollary 2.2 follows from Proposition 2.1.

\textbf{2.3. Lemma; \cite[Lemma 2.15]{Tug20}.} Let $A$ be a ring and let $M$ be a non-zero right $A$-module which is not an essential extension of its singular submodule $\text{Sing }X$. Then $M$ contains a submodule which is isomorphic to a non-zero right ideal of the ring $A$.

\textbf{2.4. Lemma.} Let $A$ be a prime right Goldie ring. 
\begin{enumerate}
\item[{\bf 1.}] 
Every non-zero ideal $B$ of the ring $A$ is an essential right ideal and contains a regular element. 
\item[{\bf 2.}] 
There exists a positive integer $n$ such that for any non-zero elements $b_1,\dots\,b_n$ of the ring $A$, the module $b_1A\oplus\dots\oplus b_nA$ contains an isomorphic copy of the free cyclic module $A_A$.
\item[{\bf 3.}] 
If the right $A$-module $X$ contains a non-singular submodule $Y$ of infinite Goldie dimension $\tau$, then $X$ contains a non-zero free submodule $F$ of infinite rank $\tau|$ and, consequently, there exists an epimorphism from the module $F$ onto any $\tau$-generated right $A$-module.
\item[{\bf 4.}] 
If there exists a non-singular is seminoetherian right $A$-module of infinite Goldie dimension $\tau$, then every $\tau$-generated right $A$-module is seminoetherian. In particular, every countably generated right $A$-module is seminoetherian.
\end{enumerate}

\textbf{Proof.} \textbf{1.} It follows from Remark 1.4 that it is sufficient to prove that $B$ is an essential right ideal. Let $C$ be a right ideal of the ring $A$ with $B\cap C=0$. Then $(AC)B\subseteq B\cap C=0$. Since the ring $A$ is prime, $C\subseteq AC=0$ and $B$ is an essential right ideal.

\textbf{2.} It follows from Remark 1.4 that the prime right Goldie ring $A$ has the classical right ring of fractions $Q$ which for some positive integer $n$ is isomorphic to ring of all $n\times n$ matrices over a division ring. We denote by $B$ the module $b_1A\oplus\dots\oplus b_n$. Then the $Q$-module $b_1Q\oplus\dots\oplus b_nQ$ contains an isomorphic copy of the module $Q_Q$. Therefore, the module $B$ contains an isomorphic copy of the module $A_A$.

{\bf 3.} It follows from Lemma 2.3 that $X$ contains a submodule which is isomorphic to a non-zero principal right ideal $bA$ of the ring $A$. It follows from \textbf{2} there exists a positive integer $n$ such that the direct sum of $n$ isomorphic copies of the module $bA_A$ contains an isomorphic copy of the free cyclic module $A_A$. Then $X$ contains a non-zero free submodule of infinite rank $|J|$.

{\bf 4.} The assertion follows from \textbf{3} and the definition of the seminoetherian module.~\hfill$\square$

\textbf{2.5. Lemma.} If a uniserial module $M$ has a maximal submodule $N$, then the module $M$ is cyclic.

\textbf{Proof.} Since the simple module $M/N$ is cyclic, there exists a cyclic submodule $X$ of $M$ which is not contained in $N$. Since $M$ is a uniserial module, $N\subseteq X$. In addition, the factor module $M/N$ is simple. Therefore, $M$ coincides with the cyclic module $X$.~\hfill$\square$

\textbf{2.6. Proposition.} Let $A$ be a prime right Goldie ring such that there exists a countably generated uniserial non-cyclic right $A$-module $E$. For a right $A$-module $M$, the following conditions are equivalent.
\begin{enumerate}
\item[{\bf 1)}] 
$M$ is a seminoetherian module;
\item[{\bf 2)}] 
$\text{Sing }M$ is a seminoetherian module and $M/\text{Sing }M$ is a finite-dimensional non-singular seminoetherian module. 
\end{enumerate}

\textbf{Proof.} We set $T=\text{Sing }M$.

2)\,$\Rightarrow$\,1). The assertion follows from Proposition 2.1.

1)\,$\Rightarrow$\,2). It follows from the definition of the seminoetherian module that the module $M/T$ is seminoetherian. By Remark 1.4, the module $M/T$ is non-singular. Let's assume that the non-singular seminoetherian module $M/T$ is infinite-dimensional. By Lemma 2.4(4), all countably generated right $A$-modules are seminoetherian. By assumption, there exists a countably generated uniserial non-cyclic right $A$-module $E$. The non-zero seminoetherian module $E$ has a maximal submodule. By Lemma 2.5, the module $E$ is cyclic; this is a contradiction.~\hfill$\square$

\textbf{2.7. Remark.} A ring $A$ is said to be \textbf{semiprimary} if its Jacobson radical is nilpotent and the factor ring $A/J(A)$ is artinian. It is directly verified that every module over a semiprimary ring is seminoetherian.

\textbf{2.8. Remark.} Let $A$ be a ring and let $M$ be a right $A$-module. It is directly verified that the module $M_A$ is seminoetherian if and only if the natural module $M_{A/B}$ is seminoetherian. Consequently, if the factor ring $A/r(M)$ is semiprimary, then it follows from Remark 2.7 that the module $M$ is seminoetherian.

\section{Modules over HNP Rings}\label{sec3}

\textbf{3.1. Proposition.} For a ring $A$, the following conditions are equivalent. 
\begin{enumerate}
\item[{\bf 1)}] 
The ring $A$ is right noetherian.
\item[{\bf 2)}] 
Every injective right $A$-module is a direct sum of uniform injective modules. 
\item[{\bf 3)}] 
All direct sums of injective right $A$-modules are injective.
\end{enumerate}

Proposition 3.1 is well known; e.g., see \cite[4.19, 4.20]{GooW89}. 

In the next Proposition 3.2 we collect for convenience a number of statements, most of which are known.

{\bf 3.2. Proposition.} Let $A$ be a non-primitive HNP ring, $E$ be an injective indecomposable non-zero singular right $A$-module ($E$ is a $P$-primary module for some primary ideal $P$), $M$ be a singular right $A$-module, and let $\{M(P_i)\}$ be the set of all primary components of the module $M$. 
\begin{enumerate}
\item[{\bf 1.}] 
In any singular right $A$-module, all finitely generated submodules are finite direct sums of cyclic uniserial modules of finite composition length. 
\item[{\bf 2.}] 
Every non-injective singular module has non-zero cyclic uniserial direct summand of finite composition length.
\item[{\bf 3.}] 
$E$ is a uniserial non-cyclic $P$-primary module without maximal submodules, all proper submodules of the module $E$ are cyclic modules of finite composition length and form a countable chain $0=X_0\subset X_1\subset \ldots \subset X_k\subset \ldots$, where $X_k/X_{k-1}$ is a simple module for any $k$ and there exists a positive integer $n$ of the module $E$) such that $X_j/X_{j-1}\cong X_k/X_{k-1}$ $\Leftrightarrow$ $j-k$ is divided by $n$ (the integer $n$ is called the \textbf{period} of $E$). 
\item[{\bf 4.}] 
If $\overline{E}$ is any non-zero homomorphic image of the module $E$, then for any cyclic submodule $\overline{X}$ of $\overline{E}$ of composition length $\ge k+n$, there exists an epimorphism $\overline{X}\to X_k$ with non-zero kernel, where $X_k$ is an arbitrary cyclic submodule from \textbf{1}.
\item[{\bf 5.}] 
If $Y$ is an arbitrary cyclic uniserial $P$-primary module, then $Y$ is isomorphic to a subfactor of the module $E$ and $Y$ is annihilated by some power of the ideal $P$.
\item[{\bf 6.}] 
If $\oplus_{i\in I}Y_i$ is a $P$-primary module, where all $Y_i$ are cyclic uniserial modules, the composition lengths of which are limited in total, then the module is annihilated by some positive power of the ideal $P$.
\item[{\bf 7.}] 
We preserve designations of from \textbf{3}. Let $\{X_k\}_{k=1}^{\infty}$ be the set of all non-zero proper (cyclic) submodules of the $P$-primary module $E$ of period $n\in\mathbb{N}$, $\{Y_{2nk}\}_{k=1}^{\infty}$ be an arbitrary set of $P$-primary cyclic modules of composition length $2nk$. Then there exist module epimorphisms $\pi\colon\oplus_{k=1}^{\infty}Y_{2nk}\to \oplus_{k=1}^{\infty}X_k$ and $h\colon\oplus_{k=1}^{\infty}X_k\to E$ and, consequently, there exists the module epimorphism $h\pi\colon\oplus_{k=1}^{\infty}Y_{2nk}\to \to E$.
\item[{\bf 8.}] 
We preserve designations from \textbf{3}. Let $\{Z_k\}_{k=1}^{\infty}$ be an arbitrary set of $P$-primary cyclic modules of composition length $d_k$ and let the integer set $\{d_k\}$ be not bounded. Then there exists an epimorphism $\oplus_{k=1}^{\infty}T_k\to E$, where all $T_k$ are cyclic $P$-primary subfactors of the module $\{Z_k\}_{k=1}^{\infty}$. In addition, the module $\oplus_{k=1}^{\infty}Z_k$ is not seminoetherian.
\end{enumerate}

\textbf{Proof.} Assertions \textbf{1}--\textbf{6} are contained in \cite{Sin74}, \cite{Sin75}, \cite{Sin76}, \cite{Sin78}, \cite{Sin79}, \cite{SinT79} or directly follow from results of these papers.

\textbf{7.} The assertion follows from assertions \textbf{3}.

\textbf{8.} It follows from \textbf{7} that there exists an epimorphism $\oplus_{i\in I}\overline{Y_i}\to E$. Let's assume that the module $\oplus_{k=1}^{\infty}Z_k$ is seminoetherian. Since there exists an epimorphism $\oplus_{k=1}^{\infty}T_k\to E$, the module $E$ is isomorphic to subfactor of the seminoetherian module $\oplus_{k=1}^{\infty}Z_k$. Therefore, the module $E$ is seminoetherian. This is a contradiction, since the module $E$ does not have maximal submodules by \textbf{3}.~\hfill$\square$

{\bf 3.3. Proposition.} Let $A$ be a non-primitive HNP ring, $M$ be a singular right $A$-module, and let $\{M(P_i)\}$ be the set of all primary components of the module $M$. The module $M$ is injective if and only if every primary component $M(P_i)$ is an injective module and also if and only if every primary component $M(P_i)$ is a direct sum of uniserial injective modules.

Proposition 3.3 follows from Propositions 3.1, 3.2(1) and 3.2(2) .

\textbf{3.4. Lemma.} Let $A$ be a non-primitive HNP ring, $M$ be a singular right $A$-module, and let $\{M(P_i)\}$ be the set of all primary components of the module $M$. 
\begin{enumerate}
\item[{\bf 1.}] 
If every primary component $M(P_i)$ of the module $M$ is seminoetherian, then the module $M$ is seminoetherian.
\item[{\bf 2.}] 
If every primary component $M(P_i)$ of the module $M$ has the non-zero annihilator, then the module $M$ is seminoetherian.
\end{enumerate}

\textbf{Proof.} With the use of Remark 1.7, assertion \textbf{1} is directly verified.

\textbf{2.} The assertion follows from assertions \textbf{1} and Remark 2.8.~\hfill$\square$

{\bf 3.5. Proposition.} Let $A$ be a non-primitive hereditary noetherian prime ring and let $M$ be a non-zero singular $A$-module.
\begin{enumerate}
\item[{\bf 1.}] 
If $M$ is a direct sum of cyclic singular modules, then every submodule of the module $M$ is a direct sum of cyclic singular modules.
\item[{\bf 2.}] 
The module $M$ contains a submodule $X$ such that $X$ is a direct sum of uniserial singular modules and the module $M/X$ is injective.
\item[{\bf 3.}] 
If $M$ is a max module, then $M$ does not contain non-zero injective submodules. \item[{\bf 4.}] 
If $M$ is a max module, then $M$ contains a submodule $X$ such that $X$ is a direct sum of cyclic uniserial singular modules and the module $M/X$ is injective.
\item[{\bf 5.}] 
If module $M$ is seminoetherian, then $M$ is a direct sum of cyclic uniserial singular modules.
\end{enumerate}

\textbf{Proof.} \textbf{1, 2.} The assertions are partial cases of \cite[Theorem 4]{Sin75} and \cite[Theorem 1]{Sin78}, respectively.

\textbf{3.} Let's assume the contrary. Then $M$ contains a non-zero injective direct summand which is a max module. It follows from Proposition 3.3 $M$ has a non-zero uniserial injective singular direct summand $N$ which is a max module. This contradicts to Proposition 3.2(3).

{\bf 4.} It follows from \textbf{2} and \textbf{3} that $M$ contains a submodule $X=\oplus_{i\in I}X_i$ such that all $X_i$ are uniserial singular non-injective modules and the module $M/X$ is injective. It follows from Proposition 3.2(5) that all modules $X_i$ are cyclic.

{\bf 5.} It follows from {\bf 4} $M$ contains a submodule $X$ such that $X$ is a direct sum of cyclic uniserial singular modules and $M/X$ is an injective singular seminoetherian module. By applying \textbf{3} to the injective singular max module $M/X$, we obtain that $M=X$.~\hfill$\square$

{\bf 3.6. Proposition.} Let $A$ be a non-primitive hereditary noetherian prime ring and let $M$ be a non-zero singular $A$-module. The following conditions are equivalent.
\begin{enumerate}
\item[{\bf 1)}] 
$M$ is a seminoetherian module.
\item[{\bf 2)}] 
Every primary component of the module $M$ is a seminoetherian module.
\item[{\bf 3)}] 
Every primary component of the module $M$ is a direct sum of cyclic uniserial modules, the composition lengths of which are limited in total.
\item[{\bf 4)}] 
Every primary component of the module $M$ has the non-zero annihilator.
\end{enumerate}

\textbf{Proof.} With the use of Remark 1.7, the equivalence of conditions \textbf{1} and \textbf{2} is directly verified. Without loss of generality, we can assume that $M$ coincides with one of its primary components. It is sufficient to prove that conditions \textbf{2}, \textbf{3} and \textbf{4 } are equivalent.

2)\,$\Rightarrow$\,3). The assertion follows from Proposition 3.2(7).

3)\,$\Rightarrow$\,4). The assertion follows from Proposition 3.2(6).

4)\,$\Rightarrow$\,2). The assertion follows from Lemma 3.4(2).~\hfill$\square$

\textbf{3.7. Remark; \cite[Proposition 1.20]{Goo76}.} For every module, any its homomorphic image with essential kernel is singular.

\textbf{3.8. Proposition.} Let $A$ be a right hereditary ring and let $Q$ be the injective hull of the module $A_A$. 
\begin{enumerate}
\item[{\bf 1.}] 
All homomorphic images of any injective right $A$-module are injective modules.
\item[{\bf 2.}]
Every indecomposable injective right $A$-module is a homomorphic image of the module $Q_A$.
\item[{\bf 3.}]
Let $Q\ne A_A$. Then $Q/A_A$ is a non-zero singular injective module.
In addition, if the ring $A$ is right noetherian, then there exist non-zero uniform singular injective modules each of them are homomorphic images of the module $Q$.
\end{enumerate}

\textbf{Proof.} \textbf{1.} The assertion is proved in \cite[Theorem 5.4]{CarE56}. 

\textbf{2.} For example, see, \cite[Corollary 8]{FacB15}.

\textbf{3.} By \textbf{1}, the module $Q/A_A$ is injective. In addition, the module $Q/A_A$ is singular by Remark 3.7. The non-zero singular injective module $Q/A_A$ is the direct sum of non-zero uniform modules $X_i$, $i\in I$, by Proposition 3.1. It follows from \textbf{2} that all modules $X_i$ are homomorphic images of the module $Q$.~\hfill$\square$

\section{Completion of the Proof of Theorem 1.1. Examples and Supplements}\label{sec4}

\textbf{4.1. Proposition.} Let $A$ be a non-primitive hereditary noetherian prime ring, $M$ be a right $A$-module, and let $T$ be the largest singular submodule of the module $M$. The following conditions are equivalent.
\begin{enumerate}
\item[\textbf{1)}]
$M$ is a seminoetherian module. 
\item[\textbf{2)}]
$T$ and $M/T$ are seminoetherian modules.
\item[\textbf{3)}]
$T$ is a seminoetherian module and $M/T$ is a finite-dimensional non-singular seminoetherian module.
\item[\textbf{4)}]
Every primary component of the module $T$ is a direct sum of cyclic uniserial modules, the composition lengths of which are limited in total, and there exists a submodule $X$ of $M$ such that $T\subseteq X$, $X/T$ is a noetherian submodule of $M/T$, and every primary component of the module $M/X$ is a direct sum of cyclic uniserial modules, the composition lengths of which are limited in total.
\item[\textbf{5)}]
Every primary component of the module $T$ is a direct sum of cyclic uniserial modules, the composition lengths of which are limited in total, $M/T$ is a non-singular finite-dimensional module, and there exists a submodule $X$ of the module $M$ such that $T\subseteq X$, $X/T$ is a noetherian hereditary essential submodule of $M/T$ which is isomorphic to a finite direct sum of uniform right ideals of the ring $A$, and every primary component of the module $M/X$ is a direct sum of cyclic uniserial modules, the composition lengths of which are limited in total.
\end{enumerate}

\textbf{Proof.} The equivalence 1)\,$\Leftrightarrow$\,2) follows from Proposition 2.1. 

2)\,$\Rightarrow$\,3). By Proposition 2.6, it is sufficient to prove that there exists a countably generated uniserial non-cyclic right $A$-module. By Proposition 3.8, there exists a non-zero uniform singular injective right $A$-module $E$. By Proposition 3.2(3), $E$ is a countably generated uniserial non-cyclic module.

3)\,$\Rightarrow$\,4). By Proposition 3.6, every primary component of the module $T$ is a direct sum of cyclic uniserial modules, the composition lengths of which are limited in total. By condition \textbf{3}, $M/T$ is a finite-dimensional non-singular seminoetherian module. Since the module $M/T$ is finite-dimensional, $M/T$ is an essential extension of some finitely generated module $X/T$. Since $A$ is a noetherian ring, $X/T$ is a noetherian module. Since $M/T$ is an essential extension of the module $X/T$, the module $M/X$ is singular. By Proposition 3.6, every primary component of the seminoetherian module $M/X$ is a direct sum of cyclic uniserial modules, the composition lengths of which are limited in total.

4)\,$\Rightarrow$\,5). Since $M/T$ is a finite-dimensional non-singular module, $M/T$ is an essential extension of some finitely generated non-singular module $N/T$. It follows from Lemma 2.3 that $N/T$ is an essential extension of some module $X/T$ which is isomorphic to a finite direct sum of uniform right ideals of the ring $A$. Then $X/T$ is a noetherian hereditary essential submodule of $M/T$ and $M/X$ is a singular seminoetherian module. By Proposition 3.6, every primary component of the seminoetherian module $M/X$ is a direct sum of cyclic uniserial modules, the composition lengths of which are limited in total.

Implications 5)\,$\Rightarrow$\,4) and 3)\,$\Rightarrow$\,2) are directly verified.~\hfill$\square$

\textbf{4.2. Corollary.} Theorem 1.1 follows from Proposition 4.1.

\textbf{4.3. Remark.} If $M$ is a uniserial max module, then it follows from Lemma 2.5 that all submodules of the module $M$ are cyclic and, in particular, $M$ is a noetherian module. 

Let $\mathbb{Z}$ be the ring of integers, $\mathbb{Q}$ be the additive group of rational numbers.

\textbf{4.4. Remark.} A module $M$ is said to be \textbf{distributive} if $X\cap(Y+Z)=X\cap Y+X\cap Z$ for any submodules $X,Y,Z$ of $M$. If all finitely generated submodules of the module $M$ are cyclic, then $M$ is called a \textbf{Bezout module}. Every uniserial module is a distributive Bezout module. It is directly verified that all subfactors of distributive modules are distributive and $\mathbb{Q}$ is a distributive non-uniserial Bezout $\mathbb{Z}$-module.

\textbf{4.5. Example.} Let $M$ be the complete pre-image in $\mathbb{Q}$ of the socle of the factor group $\mathbb{Q}/\mathbb{Z}$. Then $M/\mathbb{Z}$ is an infinite direct sum of pairwise non-isomorphic simple $\mathbb{Z}$-modules. By Theorem 1.1, $M$ is a seminoetherian distributive $\mathbb{Z}$-module which is not noetherian (we can compare this fact with Remark 4.3). In addition, $M$ is a distributive Bezout $\mathbb{Z}$-module

\textbf{4.6. Remark.} It is known that every non-zero projective module has a maximal submodule; see, e.g., \cite[Theorem 9.6.3]{Kas82}. Therefore, every hereditary module is a max module. 

\textbf{4.7. Remark.} Let $A$ be a right hereditary ring. It is known that every projective right $A$-module is hereditary; see, e.g., \cite[Theorem 5.4]{CarE56}. Therefore, it follows from Remark 4.4 that every non-zero projective right $A$-module $P$ is a max module.

\textbf{4.8. Example.} Let $M$ be a free abelian group of infinite rank. It follows from Remark 4.5 and Theorem 1.1 that $M$ be a max $\mathbb{Z}$-module which is not seminoetherian.

\textbf{4.9. Example.} There is a countable commutative  uniserial ring $A$ which is not a max $A$-module. Let $p$ be a prime integer, $\mathbb{Z}(p)$ be the localization of the ring $\mathbb{Z}$ with respect to the prime ideal $p\mathbb{Z}$, and let $M$ be the quasi-cyclic additive $p$-group. We naturally turn $M$ into a $\mathbb{Z}(p)$-module. Let $A$ be the matrix ring consisting of matrices of the form 
$\begin{pmatrix}
f&m\\
0&f
\end{pmatrix}$, where $f\in \mathbb{Z}(p)$ and $m\in M$. We note that 
$$\begin{pmatrix}
f_1&m_1\\
0&f_1
\end{pmatrix}\begin{pmatrix}
f_2&m_2\\
0&f_2
\end{pmatrix}=\begin{pmatrix}
f_1f_2&f_1m_2+m_1f_2\\
0&f_1f_2
\end{pmatrix},
$$ 
$A$ is a countable commutative ring, $M$ is an ideal of $A$ with $M^2=0$, the factor ring $A/M$ is isomorphic to the uniserial principal ideal domain $\mathbb{Z}(p)$, and $M$  is the prime radical of $A$. It is directly verified that an ideal $B$ of  $A$ contains $M$ provided $B$ is not contained in $M$. In addition, the ring $A/M$ is uniserial. Therefore, if $B$ and $C$ are two ideals of $A$ and $B\not\subseteq M$, then $B$ and $C$ are comparable with respect to inclusion. Since the module $M_{\mathbb{Z}}$ is uniserial, any two ideals of $A$, which are contained in $M$, are comparable with respect to inclusion. Thus, the ring $A$ is uniserial and all proper submodules of $M_A$  form an infinite ascending chain with respect to inclusion. Therefore, the module $M_A$ is not finitely generated and does not have maximal submodules.  Therefore, the countable commutative  uniserial ring $A$ is not a right max $A$-module.

\textbf{4.10. Example.} Let $p$ be a prime integer and  let $F=\mathbb{Z}/p\mathbb{Z}$ be the residue field modulo $p$. Then there exists a countable $F$-algebra $A$ which is a right and left max $A$-module and is not a right or left  seminoetherian  ring. Indeed, let  $G=\mathbb{С}(p^{\infty})$ be the multiplicatively written quasi-cyclic $p$-group and let $FG$ be the group algebra with augmentation ideal $M$. It is well known that $FG$ is a countable commutative local $F$-algebra with Jacobson radical $M$. Since all proper subgroups of the abelian group $G$ form an infinite ascending chain with respect to inclusion, it is well known that all proper submodules of the module $M_A$ form an infinite ascending chain with respect to inclusion, and the module $M_A$ does not have maximal submodules. Therefore, $FG]$ is a countable commutative uniserial $F$-algebra which is not a max $FG$-module. Let $A=F\langle X_1,X_2,\dots\rangle$ be the free $F$-algebra in countable set of free variables $X_1,X_2,\dots$. Then $FG$ is a homomorphic image of the algebra $A$. Since $FG$ is not a right or left seminoetherian  ring,  $A$ is not a right or left seminoetherian ring. In addition, all right (resp., left)  ideals of the algebra are free right (resp., left)  modules, see \cite{Coh64}. In particular, $A$ is a right and left hereditary ring. By Remark 4.6, the algebra $A$ is a right and left max $A$-module.

\end{document}